\documentclass{amsart}
\usepackage{amsmath,amsfonts,amssymb,amsthm,epsfig,color,tikz}
\usepackage[margin=1.25in]{geometry}

\newcommand{\bb}{\mathbb}

\newcommand{\Z}{\bb Z}
\newcommand{\R}{\bb R}
\newcommand{\N}{\bb N}

\newtheorem*{lemma*}{Lemma}
\newtheorem*{question*}{Question}

\newtheorem*{theorem*}{Theorem}

\numberwithin{Def}{section}

\newcommand{\La}{\Lambda}
\newtheorem{Theorem}{Theorem}

\newtheorem{lemma}[Theorem]{Lemma}
\numberwithin{Theorem}{section}

\begin{document}
\title[Values of Random Polynomials]{Values of Random Polynomials at Integer Points}
\author{Jayadev~S.~Athreya}
\author{Gregory A.~Margulis}
\email{jathreya@uw.edu}
\email{gregory.margulis@yale.edu}
\address{Department of Mathematics, University of Washington, Padelford Hall, Seattle, WA 98195}
\address{Department of Mathematics, Yale University, New Haven, CT, 06520}
\thanks{J.S.A partially supported by NSF CAREER grant DMS 1559860, NSF grant
   DMS 1069153, and grants DMS 1107452, 1107263, 1107367  ``RNMS: GEometric structures And Representation varieties" (the GEAR Network)."  G.A.M. is supported by NSF grant DMS 1265695. This material is based upon work while both authors were in residence at the Mathematical Sciences Research Institute in Berkeley, California during the Spring 2015 semester, supported by the National Science Foundation Grant DMS 0932078 000.}

   \begin{abstract} Using classical results of Rogers~\cite[Theorem 1]{Rogers} bounding the $L^2$-norm of Siegel transforms, we give bounds on the heights of approximate integral solutions of quadratic equations and error terms in the quantiative Oppenheim theorem of Eskin-Margulis-Mozes~\cite{EMM} for almost every quadratic form. Further applications yield quantitative information on the distribution of values of random polynomials at integral points.
   \end{abstract}
\maketitle

\section{Introduction}\label{sec:intro} \noindent Let $Q$ be an indefinite quadratic form in $n \geq 3$ variables, and suppose further that it is \emph{irrational}, that is, it is not a real multiple of a form with rational coefficients. In 1929, Alexander Oppenheim~\cite{Oppenheim} conjectured that for such forms, $Q(\Z^n)$ contains values arbitrarily close to zero.

\subsection{Small values} Our first main result bounds the size of solutions to the inequality $$\left|Q(x)\right|< \epsilon, x \in \Z^n\backslash\{0\}$$ for almost every quadratic form $Q$.

\begin{Theorem}\label{theorem:height} For every $\delta >0$, for almost every (with respect to the Lebesgue measure class) quadratic form $Q$ of signature $(p,q)$, there are constants $c_Q, \epsilon_0>0$ so that for all $\epsilon < \epsilon_0$, there is a nonzero $x \in \Z^n$ with $$\|x\| \le c_Q \epsilon^{-\left(\frac{1}{n-2} + \delta\right)}$$ and $$|Q(x)|< \epsilon.$$

\end{Theorem} 

\medskip

\subsection{Related results}\label{sec:spectral} The exponent $\frac{1}{n-2}$ is optimal. Recently, Ghosh-Kelmer~\cite{GK} (using spectral methods) and  Bourgain~\cite{Bourgain}(using analytic number theory) obtained similar (also optimal) results for ternary forms (in Bourgain's setting, diagonal ternary forms). Ghosh-Gorodnik-Nevo~\cite{GGN} have very general results on these types of problems using spectral methods. Previously, VanDerKam~\cite{VanDerKam} obtained results on the distribution of values of general homogeneous polynomials. The novelty of our results includes the elementary nature of the methods, using only classical results from the geometry of numbers, in particular, mean and variance estimates of certain counting functions on the space of lattices.

\subsection{Counting solutions} Oppenheim later strengthened his conjecture to say that $\overline{Q(\Z^n)} = \R$. This was proved by the second named author in~\cite{Margulis} in 1986, using methods from homogeneous dynamics. Subsequently, Eskin-Margulis-Mozes~\cite{EMM, EMM2} gave \emph{quantitative} versions of these results, showing, for example, that for quadratic forms $Q$ (of signature other than $(2, 1)$ or $(2, 2)$), $-\infty < a \le b <\infty$, there is a constant $c_Q$ so that $$N(Q, a, b, T) : = \# \left(Q^{-1}(a,b) \cap \Z^n \cap B(0, T) \right)$$ satisfies $$N(Q, a, b, T) \sim c_Q (b-a) T^{n-2},$$ where $\sim$ denotes that the ratio goes to $1$ as $T \rightarrow \infty$. $c_Q$ is defined by $$\left| Q^{-1}(a,b) \cap B(0, T) \right| \sim c_Q (b-a) T^{n-2},$$ where here and below $\left| \cdot \right|$ denotes the Lebesgue measure of a set in $\R^n$.We have bounds on the error term in this counting problem for almost every form $Q$.

\begin{Theorem}\label{theorem:error} For every $\delta >0$, for almost every (with respect to the Lebesgue measure class) quadratic form $Q$ of signature $(p,q)$, for any $a<b$, $$N(Q, a, b, T) = c_Q(b-a)T^{n-2} + o\left(T^{\left(\frac {n-1}{2}\right) + \delta}\right).$$

\end{Theorem}


\subsection{Lattice Point Counting}\label{sec:axiom}
\noindent These results are consequences of classical work of Rogers~\cite{Rogers} yielding variance estimates for lattice point counting as we integrate over the space of unimodular lattices, $X_n = SL(n, \R)/SL(n,\Z)$ with respect to the Haar probability measure $\mu = \mu_n$.  Another consequence of these results are general lattice point counting results. Given a subset $F \subset \R^n$ and $\La \in X_n$, let $$N(\La, F) = \#(\La \cap F).$$

\begin{Theorem}\label{dilates} Let $n \geq 4$, and $A \subset \R^n$ be measurable with $|A| = 1$. For all $\delta >0$, for almost every $\La \in X_n$, $$| N(\La, tA) - t^n | < t^{\frac 2 3 n + \delta} \mbox{ for } t >>1.$$

\end{Theorem}

\begin{Theorem}\label{sequences} Let $n \geq 3$, $B_k \subset \R^n$ be a sequence of sets, and $f(k)$ such that $\sum_{k \in \N} f(k)^{-2} < \infty.$ Then for almost every $\La \in X_n$, $$| N(\La, B_k) - |B_k|| < |B_k|^{1/2} f(k) \mbox{ for } k>>1.$$

\end{Theorem}

\medskip

\noindent\textbf{Remark:} W.~Schmidt~\cite{Schmidt} obtained similar counting results for almost all general (not unimodular) lattices. In dimension $n > 2$,  given mild technical conditions on a family of sets $\Phi$, he showed that for any $S \in \Phi$, with $|S| = V$, that for almost every $A \in GL(n, \R)$, $$|N(A\Z^n, S) - V| \le V^{1/2} \log V \psi^{1/2}(\log V),$$ where $\psi$ is a positive, non-decreasing function with $\int_{0}^{\infty} \psi^{-1}(s) ds < \infty$. He also obtained similar results for $n=2$, with $\log V$ replaced by $\log^2 V$. The distinction between our results and Schmidt's is that we restrict to \emph{unimodular lattices}, a set of measure $0$ from Schmidt's perspective. We thank the anonymous referee for pointing out that Schmidt's results give better bounds when the geometry of the set is restricted, but that our results hold for more general sets. For example, if $A$ is star-shaped, Schmidt's result (together with Fubini's theorem) gives a bound of $t^{\frac n 2 + \delta}$.

\section{Random Lattices}\label{sec:random} The key lemmas for our results are estimates of certain integrals on the space of lattices, in particular classical results of Siegel~\cite{Siegel} and Rogers~\cite{Rogers}, which estimate the \emph{mean} and \emph{variance} of the number of points of a random lattice in a set $A \subset \R^n$. 

\subsection{Siegel Transforms} We recall that given $f \in C_c(\R^n)$, its \emph{Siegel transform} is given by $\widehat{f}: X_n \rightarrow \R$ with $$\widehat{f} (\La) = \sum_{0 \neq x \in \La} f(x).$$ It is a theorem of Siegel~\cite{Siegel} that this transform is an $L^1$-isometry, that is, \begin{equation}\label{siegel}\int_{X_n} \widehat{f} d\mu = \int_{\R^n} f dm,\end{equation} where $m$ is Lebesgue measure on $\R^n$. If $f$ is the indicator function of a set $A$, this can be viewed as saying the \emph{expectation} of the number of unimodular lattice points in $A$ is simply given by $|A|$. 

\subsection{Rogers' Lemma} To compute the \emph{variance}, we have the following $L^2$ formula of Rogers which we specialize to the case where $f = I_A$, the indicator function of a measurable set $A$. Let $a := |A|$ be the Lebesgue measure of $A$, and let $B_a$ be the \emph{spherical symmetrization} of $A$, that is, the Euclidean ball around the origin $0$ with $|B_a| = |A| = a$.

\begin{lemma}\label{rogers}(\cite[Theorem 1 and 2]{Rogers}, \cite[\S4.1]{AMarg}) For $n \geq 3$, $$\int_{X_n} \widehat{I_A}^2 d\mu \le \int_{X_n} \widehat{I_{B_a}}^2 d\mu = a^2 + C_n a ,$$ where $C_n$ is as in Theorem~\ref{theorem:mink}. \end{lemma}

\noindent In particular, the \emph{variance} of the random variable $\widehat{I_A}$ is bounded by $C_n a$. 

\subsection{Random Minkowski} A consequence of Roger's lemma is the following bound on the measure of the set of unimodular lattices in $\R^n$ with large `holes', that is, the set of lattices which miss subsets $A \subset \R^n$ of large measure.  We have: 

\begin{Theorem}\label{theorem:mink}\cite[Theorem 2.2]{AMarg} Let $n \geq 3$, and $A \subset \R^n$ be a measurable set with $|A| >0$. Then $$\mu(\La: \La \cap A = \emptyset ) \le \frac{C_n}{|A|},$$ where $|A|$ denotes the Lebesgue measure of $A$ and $C_n = 8 \zeta(n-1)/\zeta(n)$.
\end{Theorem}
\medskip
\noindent Rogers' lemma combined with Chebyshev's inequality yields that for any $M >0$, \begin{equation}\label{concentration}\mu(\La: |N(\La, A) - |A||>M|A|^{1/2}  ) \le C_n M^{-2}.\end{equation}

\subsection{Proof of Theorem~\ref{sequences}} By \eqref{concentration}, $$\mu(\La: |N(\La, B_k) - |B_k||>f(k)|B_k|^{1/2}  ) \le C_n f(k)^{-2}.$$ Since $\sum_{k \geq 0} f(k)^{-2} < \infty,$ the easy part of the Borel-Cantelli lemma implies that almost every lattice $\La$ is contained in at most finitely many of the sets $$\left\{\La \in X_n: |N(\La, B_k) - |B_k||>f(k)|B_k|^{1/2}  \right\}.$$\qed

\subsection{Proof of Theorem~\ref{dilates}} Fix $\gamma  \in (1, n/3)$. We have $$\mu(\La: |N(\La, tA) - t^n| > t^{n-\gamma}) < C_n t^{-(n-2\gamma)}.$$ Let $\epsilon = t^{1-\gamma}$. Note that $$(t+ \epsilon)^n - (t- \epsilon)^n \cong 2n\epsilon t^{n-1} = 2n t^{n-\gamma}.$$  Here $\cong$ means the ratio goes to $1$ as $t \rightarrow \infty$. Between the values $t = k$ and $t=k+1$ we interpolate $k + i\epsilon, 0 < i < k^{\gamma-1}$. For each of these $k + i\epsilon$, the set we want to avoid is at most of size $k^{-(n-2\gamma)}$, and so the total measure of the set we must avoid is of size at most $$k^{-(n-2\gamma)}k^{\gamma-1} = k^{-n + 3\gamma-1}.$$ Since $\gamma < n/3$, this is summable, so proceeding as above via Borel-Cantelli, we have our result.\qed

\subsection{Proof of Theorem~\ref{theorem:height}} Let $Q_0$ be a quadratic form so that there is a $c_0 >0$ such that for all $a, b \in \R$, \begin{equation}\label{quad:volume}\left| Q_0^{-1}(a,b) \cap B(0, T) \right| = c_0 (b-a) T^{n-2} + o(T^{n-2}).\end{equation}This is a full measure condition on the set of quadratic forms (see, for example~\cite[Lemma 3.8]{EMM}). Consider the family $Q= Q_0 \circ g$, $g \in SL(n,\R)$. Note that $Q(\Z^n)$ only depends on $gSL(n,\Z)$. We first prove the theorem for almost every $gSL(n, \Z)$. Let $$A_j : =\left \{ x \in \R^n: 0<Q_0(x) < 2^{-j} , \|x\| < 2^{\left(\frac{1}{n-2} + \delta\right)j}\right\}.$$ Then $$|A_j| \cong c_0 2^{-j} \left(2^{\left(\frac{1}{n-2} + \delta\right)j}\right)^{n-2} =  c_0 2^{(n-2)\delta j}$$ Then by Theorem~\ref{theorem:mink}, $$\mu(\La: \La \cap A_j = \emptyset) < \frac{C_n}{c_0} 2^{-(n-2)\delta j}.$$ Ignoring the constants, we note that $$\sum_{j=1}^{\infty} 2^{-(n-2)\delta j} < \infty \mbox{ for all } \delta >0,$$ so for almost every $\La = g\Z^n$, there is a $j_0 = j_0(\La)$ so that for all $j \geq j_0$, $$\La \cap A_j \neq \emptyset.$$ Let $\epsilon_0 = 2^{-j_0}$. Then for all $\epsilon < \epsilon_0$, there is a nonzero $x \in \Z^n$ so that $$\|gx\| < \epsilon^{-\left(\frac{1}{n-2} + \delta\right)} \mbox{ and } Q_0(gx) < \epsilon.$$ There is a $c_g$ so that $$\|gx\| < \epsilon^{-\left(\frac{1}{n-2} + \delta\right)} \Rightarrow \|x\| < c_g\epsilon^{-\left(\frac{1}{n-2} + \delta\right)},$$ so we get the existence of nonzero $x \in \Z^n$ so that $$\|x\| < c_g\epsilon^{-\left(\frac{1}{n-2} + \delta\right)} \mbox{ and } Q(x) < \epsilon.$$ Since the volume condition (\ref{quad:volume}) holds for almost all quadratic forms $Q_0$, and we have shown that our result holds for $Q = Q_0 \circ g$ for almost every $g$, we have our result.\qed\medskip

\subsection{Proof of Theorem~\ref{theorem:error}} As above, let $Q_0$ be a quadratic form satisfying (\ref{quad:volume}) and consider the family $Q= Q_0 \circ g$, $g \in SL(n,\R)$. Let $A_0(T) = Q_0^{-1}(a,b) \cap B(0, T).$ Let $\gamma >0$. By \eqref{concentration}, $$\mu\left(\La: \left| N(\La, A_0(T)) - |A_0(T)|\right| > T^{\frac{(n-2)}{2} + \gamma}\right) < C_n \frac{\left|A_0(T)\right|}{T^{(n-2) + 2\gamma}} \cong T^{-2\gamma}.$$

\subsubsection{Transfer} Note that $$\|x\| \le T, Q(x) \in (a,b), 0 \neq x \in \Z^n  \Leftrightarrow \|g^{-1}w\| \le T, Q_0(w) \in (a,b), 0 \neq w \in g\Z^n,$$ where $w = gx$.  Let $\epsilon>0$. For $\epsilon$ sufficiently small, $$\max\left( \left|\log \|g\|_{op}\right|, \left|\log \|g^{-1}\|_{op}\right| \right) < \epsilon \Rightarrow 1-2\epsilon < \|g\|_{op}, \|g^{-1}\|_{op} < 1 + 2\epsilon.$$ So we have that for $T >>1$, $$T-1 < \|g^{-1} w\|< T \Rightarrow (T-1)(1-2\epsilon) < \|w\| < T(1+2\epsilon).$$ We will cover the space of lattices with translates of the neighborhood of identity $$\max\left( \left|\log \|g\|_{op}\right|, \left|\log \|g^{-1}\|_{op}\right|\right) < \epsilon,$$ and estimate the measure of lattices with bad error term in each set, which will then transfer to give an upper bound on the measure of the set of quadratic forms with bad error term. 

Let $\epsilon = T^{-\alpha}$. The proportion of $g \in SL(n, \R)$ with $\max\left(\left| \log \|g\|_{op} \right|, \left|\log \|g^{-1}\|_{op}\right|\right)  < \epsilon = T^{-\alpha}$ so that $\La = g\Z^n$ satisfies $$\left| N(\La, A_0(T)) - |A_0(T)|\right| > T^{\frac{(n-2)}{2} + \gamma}$$ is at most (up to constants) $$\frac {T^{-2\gamma}}{\epsilon^{n^2-1}} = T^{\alpha(n^2-1)-2\gamma}.$$ Here, the $\epsilon^{n^2-1}$ is the volume of the ball $\max\left(\left| \log \|g\|_{op}\right|, \left|\log \|g^{-1}\|_{op}\right|\right)  < \epsilon$ in the $n^2-1$ dimensional space $SL(n, \R)$. Arguing as above, there are $k^{\alpha}$ steps (intervals of size $k^{-\alpha}$) of size $\epsilon$ to go from $T=k$ to $T=k+1$, and so an upper bound on the size of the set we wish to avoid is $$k^{\alpha} k^{\alpha(n^2-1) -2\gamma} = k^{n^2\alpha - 2\gamma}.$$ For this to be summable in $k$, we require $$\gamma > \frac{1}{2} + \frac{n^2\alpha}{2}$$ So for any $\gamma > \frac{1}{2}$, by picking $\alpha$ appropriately and applying the Borel-Cantelli lemma, we get that for a set of full measure in the space of lattices, the error term is controlled by $$T^{\frac{n-1}{2} + \delta} \mbox{ where }\delta = \gamma-\frac 1 2>0.$$ Arguing again as at the end of the proof of Theorem~\ref{theorem:height}, we have our result. \qed

\section{General polynomials and volume estimates}\label{sec:volume} More generally, we consider the values of random $SL(n, \R)$-deformations of homogeneous polynomials at integer points, or, equivalently, the values of fixed homogeneous polynomials at the points of a random unimodular lattice.

We have analogs of Theorem~\ref{theorem:height} and Theorem~\ref{theorem:error} in this situation, which depend on volume estimates analogous to (\ref{quad:volume}). Given a general (homogeneous) polynomial $P$ in $n$-variables, we need to understand the $R \rightarrow \infty$ asymptotics of the volume of the region $$P^{-1}(a,b) \cap B(0, R)$$ in $\R^n$. The heuristic, if $P$ is of degree $k$, obtained via rescaling and Taylor expansion of $P$ near the $\{P=0\}$ hypersurface is that this should behave like $c_P (b-a)R^{n-k}$, where $$c_P = \int_{\{P=0\} \cap B(0, 1)} \frac{1}{\|\nabla P\|} dA,$$ where $dA$ is Lebesgue measure on the hypersurface $\{P=0\}$. However, this heuristic depends on $\nabla P$ not vanishing on $\{P=0\}$ and this variety being nonsingular, which do not hold for every possible $P$.  What we conjecture to be true in full generality is that there are constants $c'_P$, $n_P$, and $m_P$ so that the main term in this volume asymptotics is of the form $$c'_P R^{n_P} (\log(R))^{m_P}.$$ It seems likely that this could be deduced from work of Bernstein~\cite{Bernstein}; or Chambert-Loir and Tschinkel~\cite{CT} (which focuses on the asymptotics of growth of intersection of level sets (as opposed to sublevel) with balls).

%
%
\subsection{Affine Lattices}\label{sec:affine} Our results also hold in the setting of \emph{affine lattices}, with (\ref{siegel}) and Lemma~\ref{rogers} replaced by~\cite[Lemmas 3 and 4]{Athreya} respectively.

\bigskip

\noindent\textbf{Acknowledgments.}  We thank Joseph Bernstein, Peter Sarnak, and Yuri Tschinkel for useful discussions; and the anonymous referee for a careful reading of this paper. We thank the Isaac Newton Institute (INI), Cambridge, United Kingdom, where this work was initiated. In particular, we thank the organizers A.~Ghosh, A.~Gorodnik and B.~Weiss of the program on ``Interactions between Dynamics of Group Actions and Number Theory", June/July 2014 at INI, as well as the Mathematical Sciences Research Institute for its hospitality in Spring 2015 as part of the program on ``Arithmetic and Geometric aspects of Homogeneous Dynamics". J.S.A. thanks the organizers of the Stellenbosch/African Institute of Mathematical Sciences conference on Number Theory, January 2015 and Yale University for its hospitality. G.A.M. thanks the University of Illinois for its hospitality.
\bigskip


\end{document}